\documentclass[a4paper,12pt]{amsart}
\usepackage[T1]{fontenc}
\usepackage[utf8]{inputenc}
\usepackage{lmodern}
\usepackage[english]{babel}
\usepackage{hyperref}
\usepackage{amsmath, amsthm}
\usepackage{color}

\newtheorem{conj}{Conjecture}
\newtheorem{prop}{Proposition}

\newtheorem{theo}{Theorem}

\newcommand{\R}{\mathbf R}
\newcommand{\C}{\mathbf C}

\newcommand{\PGL}{\mathrm{PGL}}

\newcommand{\geom}{\mathrm{geom}}
\newcommand{\hyp}{\mathrm{hyp}}
\newcommand{\periph}{\mathrm{periph}}
\newcommand{\Hol}{\mathrm{Hol}}
\newcommand{\Vol}{\mathrm{Vol}}

\title[Volume of representations]{Volume of representations and birationality of peripheral holonomy}
\author{Antonin Guilloux}

\begin{document}

\begin{abstract}
We discuss here a generalization of a theorem by Dunfield stating that the peripheral holonomy map,
from the character variety of a 3-manifold to the A-polynomial is birational. 
Dunfield's proof involves the rigidity of maximal volume. The volume is still an important ingredient in this paper.
Unfortunately at this point no complete proof is done.
 Instead, a conjecture is stated about the volume function on the character variety 
 that would imply the generalized birationality result.
\end{abstract}

\maketitle

\section*{Introduction}

Let $M$ be an orientable hyperbolic manifold with one cusp (e.g. a knot complement) and $\Gamma$ its fundamental group. Choose an embedding $\mathbf Z^2\to \Gamma$ of the fundamental group of the peripheral torus. Let $X(\Gamma,\PGL(2,\C))$ be the character variety:
$$ X(\Gamma,\PGL(2,\C)) = \mathrm{Hom}(\Gamma,\PGL(2,\C))// \PGL(2,\C) $$
and $X_2$ be the component of the hyperbolic monodromy $\rho_\mathrm{hyp} : \Gamma\to \PGL(2,\C)$. Moreover, let $\Hol_\mathrm{periph}$ denote, as in \cite{Guillouxperiph,FalbelGuillouxKoseleffRouillier}, the peripheral holonomy. It is the map $X(\Gamma,\PGL(2,\C))\to X(\mathbf Z^2,\PGL(2,\C))$ naturally induced by the restriction of a representation $\rho:\Gamma\to \PGL(2,\C)$ to the peripheral $\mathbf Z^2$. Note that the image of $X_2$ by the map $\Hol_\periph$ is the usually computed $A$-polynomial when $M$ is a knot complement.
In his paper \cite{Dunfield}, Dunfield proves the following theorem:
\begin{theo}\label{thm:Dunfield}
The map $\Hol_\mathrm{periph}$, from $X_2$ to its image, is a birational map.
\end{theo}

We discuss in this paper a possible generalization of this result
 to the case of target group $\PGL(n,\C)$ for $n\geq 2$ and multicusped $M$.
  The generalization for $n=2$ and $M$ multicusped is proven by Klaff-Tillmann 
  \cite{KlaffTillman}.

Dunfield's proof uses in a crucial way the properties of the volume of representations. We will review the needed facts about this function. A major obstacle for a generalization is that the proof also uses -- as does \cite{KlaffTillman} -- the fact that the set of points in $X_2$ corresponding to the hyperbolic monodromy of a Dehn surgery of $M$ is Zariski-dense in $X_2$. Rigidity of volume maximality for these points then grants the theorem. Our main problem is that this Zariski density does not hold for $n>2$.

We present in this paper a possible workaround, still using properties of the volume. Indeed we will state a conjecture, of geometric flavour,  that implies the birationality in Dunfield's theorem. Hopefully, the conjecture may be easier to tackle in the general case than the birationality problem.

We will recall in section \ref{s:charvar} some facts about the character variety and the map $\Hol_\periph$ following mainly the presentation of \cite{FalbelGuillouxKoseleffRouillier}. One important result for the present paper is the local rigidity theorem (see theorem \ref{thm:localrigidity}). We define the geometric representation as the composition $\rho_\geom= r_n\circ \rho_\hyp$, where $r_n$ is the unique irreducible representation from $\PGL(2,\C)$ to $\PGL(n,\C)$. The component $X_n$ of the character variety $X(\Gamma,\PGL(n,\C))$ containing the geometric representation is called the \emph{geometric component}.

 We proceed in section \ref{s:volume} with the definition of the \emph{volume} of representations using bounded cohomology, following Bucher-Burger-Iozzi \cite{BucherBurgerIozzi}. Moreover we recall another approach to this function: the more combinatorial notion of volume defined 
 originally by Thurston for $n=2$, then by Bergeron-Falbel-Guilloux 
 \cite{BergeronFalbelGuilloux} for $n=3$ and fully generalized both by Dimofte-Gabella-Goncharov and Garoufalidis-Goerner-Zickert \cite{DimofteGabellaGoncharov, GaroufalidisGoernerZickert}.
 
Both approaches to the volume function yield valuable informations: first, through bounded cohomology one gets the volume rigidity of the geometric representation. Indeed, as proven in \cite{BucherBurgerIozzi}, the volume map has a unique maximum on the geometric component $X_n$, which corresponds to $\rho_\geom$.
An information given by the combinatorial approach to the volume is a formula for its derivative which is only expressed in terms of the peripheral representations (see theorem \ref{thm:Volderivative}). In other words, the volume map on $X_n$ -- as in the case $n=2$ -- always factors through $\Hol_\periph$ (see proposition \ref{prop:volumefactors}) on a Zariski-open set.

A problem that seems rather natural is to explicitly study the volume as a function defined on $X_n$ and try to retrieve geometrical information from its behaviour. 
For example, it raises the following question, which will be relevant for us is the following:  
if a representation has almost the maximal volume, is it almost the geometric 
representation? We state the following conjecture:
\begin{conj}\label{question}
Let $M$ be an orientable cusped hyperbolic $3$-manifold, $X_n$ the geometric component of the character variety with target group $\PGL(n,\C)$ and $\Hol_\periph$ the peripheral holonomy map. Then, outside of a neighbourhood of $[\rho_\geom]$, the volume is bounded away from its maximum on $X_n$.
\end{conj}

This conjecture raises an interesting and natural question \emph{per se}. But it also turns out that this conjecture leads to a generalization of Dunfield's theorem. Indeed, we prove the following result (see section \ref{s:proof}):

\begin{theo}\label{thm:equiv}
Let $M$ be an orientable cusped hyperbolic manifold, $X_n$ the geometric component of the character variety with target group $\PGL(n,\C)$ and $\Hol_\periph$ the peripheral holonomy map. 
Suppose that outside of a neighborhood of $[\rho_\geom]$, the volume is bounded away from its maximum on $X_n$. 

Then the map $\Hol_\periph$ is a birational isomorphism between $X_n$ and its image.
\end{theo}

Some experimental evidences for this conjecture have been gathered and
 are presented in the last section \ref{s:exp}

\section{Character variety and local rigidity}\label{s:charvar}

Most of the material of this section is already reviewed, with  the same notations, in \cite{FalbelGuillouxKoseleffRouillier}. Let $n\geq 2$ be an integer. Let $N$ be a compact, oriented $3$-manifold, with non-empty boundary $\partial N$, such that its interior $M$ is an oriented cusped hyperbolic $3$-manifold. Denote by $\Gamma$ the fundamental group of $M$.

\subsection{Character variety and geometric representation}

Let $G$ be a finitely generated group. The character variety $X(G,\PGL(n,\C))$ is the GIT quotient:
$$X(G,\PGL(n,\C)) = \mathrm{Hom}(G,\PGL(n,\C))//\PGL(n,\C).$$
We refer to a paper bu Sikora \cite{Sikora} for a general exposition of this object. We will restrict to two cases: first when $G=\Gamma=\pi_1(M)$ and second when $G=\mathbf Z^2$. We know in this setting that the character variety is an affine algebraic variety.

There is always a distinguished point in $X(\Gamma,\PGL(n,\C))$: the class $[\rho_\geom]$ of the \emph{geometric representation}. This representation is defined as the composition of the hyperbolic monodromy $\rho_\hyp : \Gamma \to \PGL(2,\C)$ of $M$ with the (unique) irreducible representation\footnote{Recall for example that when $n=3$, the representation $r_n$ is also known as the adjoint representation.} $r_n : \PGL(2,\C) \to \PGL(n,\C)$:
$$\rho_\geom = r_n \circ \rho_\hyp.$$

The character variety is not irreducible and we will not study it entirely. The main object of interest in this paper will be the \emph{geometric component} $X_n$: it is the unique component of $X(\Gamma,\PGL(n,\C))$ that contains $[\rho_\geom]$.

\subsection{Peripheral holonomy and local rigidity}

Let $t$ be the number of peripheral tori\footnote{The reader may as well assume $t=1$ and 
restricts to the case of a knot complement. It will not really interfere, and may simplify notations.} of $M$. 
Let $(T_i)_{1\leq i\leq t}$ be the collection of these tori and, for each $i$, $\Delta_i\simeq \mathbf 
Z^2$ be (a choice of) an injection of $\pi_1(T_i)$ inside $\Gamma$.

For any representation $\rho :\Gamma \to \PGL(n,\C)$, one may consider its restrictions $\rho_{\Delta_i}$ to the various subgroups $\Delta_i$. This restriction defines a natural algebraic map, called \emph{peripheral holonomy}:
$$\Hol_\periph : X(\Gamma,\PGL(n,\C)) \to \prod_{i=1}^t X(\Delta_i,\PGL(n,\C)).$$

This map is, together with the volume, a main character of this paper. Dunfield's theorem, indeed, states that (in the case $n=2$, $t=1$) it is a birational isomorphism between $X_n$ and its image. This map has already been studied through different points of view \cite{Guillouxperiph,MenalFerrerPorti,BergeronFalbelGuillouxKoseleffRouillier}.
For the scope of this paper, we need to recall the result of Menal-Ferrer and Porti \cite{MenalFerrerPorti} (see also \cite{BergeronFalbelGuillouxKoseleffRouillier, GuillouxPGLn}):
\begin{theo}[Menal-Ferrer -- Porti]\label{thm:localrigidity}
On a neighborhood of $[\rho_\geom]$ in $X_n$, the map $\Hol_\periph$ is a bijection on its image.
\end{theo}
The theorem, as proven in the references, is more precise than this statement, giving local parameter for $X_n$ around $[\rho_\geom]$. We will not need the enhanced version.

As a corollary, one may note that $\Hol_\periph$ is a ramified covering (it is of finite degree). We would like to prove that its degree is indeed $1$.

\section{Volume through bounded cohomology or combinatorics}\label{s:volume}

\subsection{Bounded cohomology and rigidity}

Recall the important definitions and results from the work of Bucher-Burger-Iozzi 
\cite{BucherBurgerIozzi} about the notion of \emph{volume} of a representation 
$\rho : \Gamma \to \PGL(n,\C)$. In the cited paper, the volume map -- there called 
Borel invariant -- is defined for any $[\rho]$ in the character variety. It is the evaluation 
on the fundamental class of $N$ in $H_3(N,\partial N)$
 of a suitably constructed bounded cocycle on 
$\PGL(n,\C)$ pulled back by $\rho$.  
We will not review here the precise definition. For the present work, the approach of 
the cited article 
gives a crucial theorem, namely the volume rigidity result: the maximal volume is only
attained once on the whole character variety, at $[\rho_\geom]$. 
This rigidity theorem is a key point for the 
present paper. Recall that $\Vol_\hyp$ is the hyperbolic volume of $M$

\begin{theo}[Bucher-Burger-Iozzi]\label{thm:BBI}
The map $$\Vol : X(\Gamma,\PGL(n,\C))\to \R$$ is onto $[-\frac{n(n^2-1)}{6}\Vol_\hyp;\frac{n(n^2-1)}{6}\Vol_\hyp]$. Moreover, for any point $[\rho]$ in  $X(\Gamma,\PGL(n,\C))$, we have
\begin{itemize}
\item $\Vol([\rho])=\frac{n(n^2-1)}{6}\Vol_\hyp$ iff $[\rho]=[\rho_\geom]$, 
\item $\Vol([\rho])=-\frac{n(n^2-1)}{6}\Vol_\hyp$ iff $[\bar\rho]=[\rho_\geom]$.
\end{itemize}
\end{theo} 

The stated conjecture \ref{question} comes from a question arisen during the work here 
presented: is it possible to "perturbate" the previous theorem. In other terms, does it 
holds that a representation has almost the hyperbolic volume if and only if it is almost 
the geometric one. Some experimentations to check
this conjecture in the case $n=2$ are presented in the last section \ref{s:exp}

\subsection{Combinatorics and computation of the derivative}

Another approach for the volume function (and historically the first one) was proposed by Thurston first, then in 
\cite{BergeronFalbelGuilloux} in the case $n=3$ and generalized for any $n$ in \cite{GaroufalidisGoernerZickert,DimofteGabellaGoncharov}. . It is 
combinatorial and goes through a triangulation of $M$. The idea is to work with representations decorated by flags. For a decorated representation, each tetrahedron of the 
triangulation becomes a tetrahedron of flags (a hyperbolic ideal tetrahedron in the case $n=2$). For these tetrahedra, a notion of volume is defined by sums of Bloch-Wigner dilogarithms of cross-ratios.

We will not describe more precisely this approach. Still, 
two valuable informations are:
\begin{enumerate}
\item The combinatorial notion of volume (defined on a Zariski-open subset of the character variety) coincide with the map $\Vol$ defined in the previous section, as explained in \cite[Section 2]{BucherBurgerIozzi}.
\item As a consequence, the map $\Vol$ is real analytic on a Zariski open subset of 
$X(\Gamma,\PGL(n,\C))$ and we know a formula for its first derivative.
\end{enumerate}
Indeed, an important achievement of the combinatorial approach is that it yields a formula for 
the derivative of the volume. And the crucial point for us is that this formula only depends on 
the peripheral holonomy. Indeed, the following theorem is proven by Neumann-Zagier 
\cite{NeumannZagier} for $n=2$, in \cite[Section 11.1]{BergeronFalbelGuilloux} for $n=3$ and its 
generalization to any $n$ is given in \cite{GuillouxPGLn}.

\begin{theo}\label{thm:Volderivative}
There is a $1$-form $\mathrm{wp}$ on $\prod_{i=1}^t X(\Delta_i,\PGL(n,\C))$ such that $d\Vol$ is 
the pullback by $\Hol_\periph$ of $\mathrm{wp}$ on a Zariski-open subset of $X_n$.
\end{theo}

\section{The birationality result: a conditional proof}\label{s:proof}

We prove in this section theorem \ref{thm:equiv}. Let $M$ be, as before, 
an oriented cusped hyperbolic manifold, $t$ the number of its cusps. 
Recall that $X_n$ is the geometric component of its character variety for $\PGL(n,\C)$. 

Throughout this section, we assume the conjecture \ref{question} holds for this particular $M$: outside a neighborhood of $[\rho_\geom]$ in $X_n$, the volume function $\Vol$ is bounded away from its maximum on $X_n$.

Under this assumption, we prove theorem \ref{thm:equiv} stating that
$\Hol_\periph$ is a birational isomorphism between $X_n$ and its image.
The first step of the proof is already proven -- and crucial -- when $n=2$ in \cite{Dunfield} and \cite{KlaffTillman}:
\begin{prop}\label{prop:volumefactors}
There is a real-analytic map $V$ from a Zariski-open subset of $\Hol_\periph(X_n)$ in $\prod_{i=1}^t X(\Delta_i,\PGL(n,\C))$ to $\R$ such that for any $[\rho]$ in a Zariski-open subset of $X_n$, $\Vol([\rho])=V(\Hol_\periph([\rho]))$.
\end{prop}

\begin{proof}
The point is to prove that $\mathrm{wp}$ is exact on a Zariski-open subset of the image of 
$X_n$. $V$ is then one primitive.

Consider a loop $l$ in $\Hol_\periph(X_n)$ and assume it avoids the ramification locus of 
$\Hol_\periph$. Let $\bar l$ be a lift in $X_n$. The two ends of $\bar l$ have a volume 
differing by the integral $\int_l \mathrm{wp}$. If $\bar l$ is not a loop, 
we may continue the lifting of $l$ to construct a sequence of points in $X_n$ such 
that two consecutive points always have volumes differing by this integral. 
As the volume is bounded on $X_n$, this forces the integral to vanish.

Hence $\mathrm{wp}$ is exact on a Zariski-open subset of the image of $X_n$ and the function $V$ is its primitive whose value at $\Hol_\periph([\rho_\geom])$ is $\frac{n(n^2-1)}{6}\Vol_\hyp$.
\end{proof}

Let us proceed with the proof of theorem 2.
\begin{proof}
Let $U$ be a neighborhood of $[\rho_\geom]$ in $X_n$ such that:
\begin{itemize}
\item Restricted to $U$, the map $\Hol_\periph$ is a bijection onto its image (local rigidity, see thm \ref{thm:localrigidity}).
\item $\Vol^{-1}(\Vol(U))=U$.
\end{itemize}
 The fact that $U$ exists is a consequence of the conjecture: $\Vol(U)$ is a small neighbourhood of $\frac{n(n^2-1)}{6}\Vol_\hyp$ and it has no preimage far from $[\rho_\geom]$.
 
 Now, let $z$ be in $\Hol_\periph(U)$ such that $V$ is defined at $z$. We want to prove that $z$ has a unique preimage by $\Hol_\periph$. Let $[\rho_0]\in U$ be such that $z=\Hol_\periph([\rho_0])$. By definition $V(z)=\Vol([\rho_0])$. Let $[\rho]$ be any point in $X_n$ such that $\Hol_\periph([\rho])=z$.
As $\Vol([\rho])=V(z)$ by the previous proposition, we get that $[\rho]$ belongs to $U$. By definition of $U$, it implies that
$[\rho]=[\rho_0]$.

Hence there is an open subset in the image $\Hol_\periph(X_n)$ on which the preimages of points are singletons. This means that the degree of $\Hol_\periph$, from $X_n$ to its image, is $1$.
\end{proof}

\section{Experimental evidences}\label{s:exp}

We have implemented, using Sage 7.1 \cite{sage} and Snappy \cite{SnapPy}, a threefold test 
to explore the validity of 
the conjecture in the case $n=2$. Snappy contains a census of 
61 911 cusped oriented manifold build as the gluing of at most 9 tetrahedra. We focus on
the  $1 263$ manifolds with at most $6$ tetrahedra.

Here is a description of the test built from three different tests
of increasing complexity and power. Let $M$ be an orientable manifold 
and $V_M$ its hyperbolic volume. We assume that $M$ is ideally triangulated by $\nu$ 
tetrahedra. Let $v_0 \simeq 1.015$ be the volume of the regular ideal hyperbolic 
tetrahedron. It is the maximal volume for an ideal hyperbolic tetrahedron.

The computation will be done
at the natural level for a Snappy object: the deformation variety, defined by the famous gluing equations of Thurston (see for example \cite{Tillmann}). We work with decorated 
representations, i.e. monodromies of
the gluing of hyperbolic realisations for the $n$ tetrahedra. 
The deformation variety is seen as an affine algebraic subset 
in $\C^\nu$ (often written in $\C^{3\nu}$ to keep tracks of the 
symmetries of the tetrahedra): each tetrahedron is described 
by the cross-ratio of its vertices.

A crucial point is that the volume is the sum of Bloch-Wigner dilogarithm of the 
cross-ratios and hence extend to the compactification of the character variety in 
$(\mathbf{CP}^1)^\nu$, as the dilogarithm is well-defined and continuous 
on $\mathbf{CP}^1$. We can check the conjecture for $M$ by proving that, at ideals points in this compactification, the generalized volume is bounded away from $V_M$.

Now define $k$ to be $\nu - \lceil\frac{V_M}{v_0}\rceil + 1$. The meaning of $k$ is: if we 
know for an (ideal) decorated representation that the volume of $k$ tetrahedra among the 
$\nu$ vanish, then the volume of this (ideal) representation is less than $V_M$: indeed, even 
putting the remaining $\nu-k$ tetrahedra to the maximal volume $v_0$, we do not reach $V_M$. 

\medskip

{\bf First test \\}
 The first test is very crude: if $k\leq 2$ then the conjecture holds for 
$M$. Indeed, at an ideal point at least a tetrahedron degenerates and has 
volume $0$. But, one may further assume that another tetrahedron is 
non-positively oriented, which implies it has of volume $\leq 0$. This assumption is 
licit because on 
the subset of the deformation variety where each tetrahedron is positively 
oriented, the 
volume function is convex in 
suitable coordinates. Hence, at the boundary of the positively oriented part, 
the volume 
of ideal points is strictly less than the maximum $V_M$. We need only to check 
the volume of ideal points outside of this boundary.

This first test grants that the conjecture holds for $1 144$ among the $1 263$ manifolds. Moreover, as it is straightforward to compute, one may check that $25 986$ out of the 
$61 911$ pass the test.

\medskip

{\bf Second test \\}
For the $119$ manifolds left undecided by the first test, we may use the 
logarithmic limit set to determine the minimal number $l$ of tetrahedra 
degenerating at an ideal point. The reader may find in \cite{Tillmann} a 
presentation of the logarithmic limit set for the deformation variety. Recall 
that any point in the logarithmic limit set corresponds to ideal points for 
the deformation variety.

Snappy may be used to recover the gluing equations defining the deformation 
variety. Then SageMath, through the software Gfan  \cite{gfan}, is able to 
explicitly compute the points in the logarithmic limit set.
Remark that any tetrahedron whose
coordinate in the logarithmic limit set does not vanish 
does indeed degenerate at 
the corresponding ideal points.
Thus, we compute the minimum $l$
of degenerating tetrahedra at an ideal point. If $l\geq k$, then the manifold 
pass the test.

It is a crude test, as each non degenerating tetrahedron is set to the 
maximal volume.
Still further 47 manifolds pass the test. Note that this test, as written, may not be executed for manifolds with more tetrahedra. Indeed for more than $7$ tetrahedra, Gfan does not compute the logarithmic limit. The system defining the deformation variety seems too big ($\geq 21$ variables and equations).

\medskip

The remaining $72$ manifolds all have $6$ tetrahedra, and have all a single cusp. We try on them a third test, more involved computationally.

\medskip

{\bf Third test \\}
We now try to compute explicitly the ideal points in $(\mathbf{CP}^1)^\nu$. This
problem is hard in general (recall that we have a system with $18$ variables, $18$ equations and of degree around $10$). A partly hand-driven computation is 
still often possible.

We compute Gr\"obner basis (with the Giac 
Gr\"obner engine \cite{giac} which appears to be the most effective for this problem), 
for the ideal
defining the deformation variety and then the ideal $I$ defining its ideal 
points in
$(\mathbf{CP}^1)^\nu$. This part of the computation is done in a spirit similar to
\cite{FalbelGuillouxKoseleffRouillier} (trying to project on few variables and then reconstructing the 
whole ideal).

When $I$ is $0$-dimensional, we may then use the rational univariate 
representation \cite{RUR}, as for example in \cite{FalbelKoseleffRouillier} 
through a call to the relevant Maple function to get a parametrization of ideal points. We are able to compute compute explicitly an approximated volume 
at ideal points by computing dilogarithms.

This procedure works for $29$ from the $72$ manifolds. And in each successful case we 
check that the maximal volume for ideal points is indeed less than the hyperbolic volume $V_M$. In fact, the volumes computed never exceed $10^{-13}$.

There are two risks of failure for this procedure: it happens that the ideal 
$I$ is not 
$0$-dimensional. It also happens that the computation is too long. In the first 
case of failure, some further study has been done on one example, enabling to 
check the conjecture for this example. But such a study is not at all automatic.

\bigskip

At the end of this threefold test, $1 221$ out of $1263$ pass the test, giving hope for the conjecture. For the 
remaining $42$ ones and manifolds with more than $7$ tetrahedra, the computations as presented here are too complicated.

\bibliographystyle{amsplain}
\bibliography{biblio}

\end{document}